\documentclass[twoside,11pt,leqno]{article}
\usepackage{amssymb}
\usepackage{amsthm}
\usepackage[tbtags]{amsmath}
\usepackage{doc}
\usepackage{latexsym}
\usepackage{amscd}
\font\headd=cmr8 
\theoremstyle{change}

\begin{document}
\thispagestyle{plain}
 \markboth{}{}
\small{\addtocounter{page}{546} \pagestyle{plain}
\noindent{\scriptsize KYUNGPOOK Math. J. 43(2003), 547-566}
\vspace{0.2in}\\
\noindent{\large\bf A Note on Maass-Jacobi Forms}
\footnote{{}\\
\indent Received October 21, 2002.
\indent This work was supported by Korea Research Foundation
Grant(KRF-2000-041-D00005).}
\vspace{0.15in}\\
\noindent{\sc Jae-Hyun Yang}
\newline
{\it Department of Mathematics, Inha University, Incheon 402-751,
Korea\\
e-mail} : {\verb|jhyang@inha.ac.kr|}
\vspace{0.15in}\\
{\footnotesize {\sc Abstract.} In this paper, we introduce the
notion of Maass-Jacobi forms and investigate some properties of
these new automorphic forms. We also characterize these
automorphic forms in several ways.}
\vspace{0.2in}\\
\noindent{\bf 1. Introduction} \setcounter{equation}{0}
\renewcommand{\theequation}{1.\arabic{equation}}
\vspace{0.1in}\\
\indent We let $SL_{2, 1}({\mathbb R})=SL(2,{\mathbb R})\ltimes
{\mathbb R}^{(1, 2)}$ be the semi-direct product of the special
linear group $SL(2, {\mathbb R})$ of degree $2$ and the
commutative group ${\mathbb R}^{(1, 2)}$ equipped with the
following multiplication law
\begin{equation}
(g,\alpha)\ast (h, \beta)=(gh,\,\alpha\,^th^{-1}+\beta),\ \ \
g,h\in SL(2,{\mathbb R}),\ \ \alpha,\beta\in {\mathbb R}^{(1, 2)},
\end{equation}
where ${\mathbb R}^{(1,2)}$ denotes the set of all $1\times 2$
real matrices. We let
$$SL_{2,1}({\mathbb Z}) = SL(2,{\mathbb Z})\ltimes {\mathbb Z}^{(1,2)}$$
be the discrete subgroup of $SL_{2,1}({\mathbb R})$ and $K=SO(2)$
the special orthogonal group of degree $2$.

\medskip

Throughout this paper, for brevity we put
$$G = SL_{2, 1}({\mathbb
R}),\ \ \Gamma_1=SL(2,{\mathbb Z})\ \ \ \text{and}\ \ \
\Gamma=SL_{2,1}({\mathbb Z}).$$ Let ${{\mathbb H}}$ be the
Poincar{\'e} upper half plane. Then $G$ acts on ${\mathbb
H}\times{\mathbb C}$ transitively by
\begin{equation}
(g,\alpha)\circ
(\tau,z)=((d\tau-c)(-b\tau+a)^{-1},(z+\alpha_1\tau+\alpha_2)
(-b\tau+a)^{-1}),
\end{equation}
where $g=\begin{pmatrix} a & b \\ c & d
\end{pmatrix}
\in SL(2,{\mathbb R}),\ \alpha= (\alpha_1,\alpha_2)\in {\mathbb
R}^{(1,2)}$ and $(\tau, z)\in {\mathbb H}\times{\mathbb C}.$ We
observe that $K$ is the stabilizer of this action (1.2) at the
origin $(i,0)$. ${\mathbb H}\times{\mathbb C}$ may be identified
with the homogeneous space $G/K$ in a natural way.

\medskip

The aim of this paper is to define the notion of Maass-Jacobi
forms generalizing that of Maass wave forms and study some
properties of these new automorphic forms. For the convenience of
the reader, we review Maass wave forms. For $s\in {\mathbb C}$, we
denote by $W_s(\Gamma_1)$ the vector space of all smooth bounded
functions $f:SL(2, {\mathbb R})\longrightarrow {\mathbb C}$
satisfying the following conditions (a) and (b)\,:
\smallskip

\pagestyle{myheadings}
 \markboth{\headd Jae-Hyun Yang $~~~~~~~~~~~~~~~~~~~~~~~~~~~~~~~~~~~~~~~~~~~~~~~~~~~~~~~~~~$}
 {\headd $~~~~~~~~~~~~~~~~~~~~~~~~~~~~~~~~~~~~~~$ A Note on Maass-Jacobi Forms}

(a)\ \ $f(\gamma gk)=f(g)$ \ for all $\gamma\in \Gamma_1,\ g\in
SL(2,{\mathbb R})$ and $k\in K.$

\smallskip

(b)\ \ $\Delta_0 f\,=\,{{1-s^2}\over 4}\,f,$

\smallskip

\noindent where $\Delta_0 =\,y^2\left(\,{{\partial^2}\over
{\partial x^2}}+ {{\partial^2}\over {\partial y^2}}\right)\,-\,
y{{\partial^2}\over {\partial x
\partial\theta}}+{\frac 54}{{\partial^2}\over {\partial\theta^2}}$
is the Laplace-Beltrami operator associated to the $SL(2,{\mathbb
R})$-invariant Riemannian metric
$$ds_0^2 = {1\over {y^2}}(dx^2 +
dy^2)+ \left( d\theta +{{dx}\over {2y}}\right)^2$$ on
$SL(2,{\mathbb R})$ whose coordinates $x, y, \theta\,(\,x\in
{\mathbb R},\,y>0,\,0\leq \theta < 2\pi\,)$ are given by
$$g=\begin{pmatrix} 1 & x\\ 0 & 1\end{pmatrix} \begin{pmatrix} y^{1/2} & 0\\ 0 &
y^{-1/2}\end{pmatrix} \begin{pmatrix} \cos\,\theta &
\sin\,\theta\\ -\sin\,\theta & \cos\,\theta\end{pmatrix},\ \ \
g\in SL(2,{\mathbb R})$$ by means of the Iwasawa decomposition of
$SL(2, {\mathbb R}).$ The elements in $W_s(\Gamma_1)$ are called
{\it Maass\ wave\ forms}. It is well known that $W_s(\Gamma_1)$ is
nontrivial for infinitely many values of $s.$ For more detail, we
refer to [6], [9], [13], [17] and [20].

\medskip

The paper is organized as follows. In Section 2, we calculate the
algebra of all invariant differential operators under the action
(1.2) of $G$ on ${\mathbb H}\times{\mathbb C}$ completely. In
addition, we provide a $G$-invariant Riemannian metric on
${\mathbb H}\times{\mathbb C}$ and compute its Laplace-Beltrami
operator. In Section 3, using the above Laplace-Beltrami operator,
we introduce a concept of Maass-Jacobi forms generalizing that of
Maass wave forms. We characterize Maass-Jacobi forms as smooth
functions on $G$ or $S{{\mathcal P}}_2\times {\mathbb R}^{(1,2)}$
satisfying a certain invariance property, where $S{{\mathcal
P}}_2$ denotes the symmetric space consisting of all $2\times 2$
positive symmetric real matrices $Y$ with $\det \,Y\,=\,1.$ In
Section 4, we find the unitary dual of $G$ and present some
properties of $G$. In Section 5, we describe the decomposition of
the Hilbert space $L^2(\Gamma\backslash G)$. In the final section,
we make some comments on the Fourier expansion of Maass-Jacobi
forms.
\vspace{0.1in}\\
\noindent {\bf Notations.} We denote by ${\mathbb Z},\ {\mathbb
R}$ and ${\mathbb C}$ the ring of integers, the field of real
numbers and the field of complex numbers respectively. ${\mathbb
Z}^+$ denotes the set of all positive integers. $F^{(k, l)}$
denotes the set of all $k \times l$ matrices with entries in a
commutative ring $F$. For a square matrix $A,\ \sigma(A)$ denotes
the trace of $A$. For any $M\in F^{(k,l)},\ ^tM$ denotes the
transpose of $M$. For $A\in F^{(k,l)}$ and $B\in F^{(k, k)},$ we
set $B[A]=\,^tABA.$ We denote the identity matrix of degree $n$ by
$E_n$. ${{\mathbb H}}$ denotes the Poincar{\' e} upper-half plane.
\vspace{0.2in}\\
\noindent{\bf 2. Invariant Differential Operators  on ${{\mathbb
H}}\times {\mathbb C}$} \setcounter{equation}{0}
\renewcommand{\theequation}{2.\arabic{equation}}
\vspace{0.1in}\\
\indent We recall that $S{{\mathcal P}}_2$ is the symmetric space
consisting of all $2\times 2$ positive symmetric real matrices $Y$
with $\det\,Y=1.$ Then $G$ acts on $S{\mathcal P}_2\times {\mathbb
R}^{(1,2)}$ transitively by \begin{equation} (g,\alpha)\cdot
(Y,V)=\,(\,gY\,^t\!g,\,(V+\alpha)\,^t\!g\,),
\end{equation}
where $g\in SL(2,{\mathbb R}),\ \alpha\in{\mathbb R}^{(1, 2)},\
Y\in S{{\mathcal P}}_2$ and $V\in{\mathbb R}^{(1, 2)}.$ It is easy
to see that $K$ is a maximal compact subgroup of $G$ stabilizing
the origin $(E_2,0).$ Thus $S{\mathcal P}_n\times {\mathbb R}^{(m,
n)}$ may be identified with the homogeneous space $G/K$ as
follows\,:
\begin{equation}
G/K\ni (g,\alpha)K\longmapsto (g,\alpha)\cdot (E_2,0)\in
S{\mathcal P}_2\times {\mathbb R}^{ (1,2)},
\end{equation}
where $g\in SL(2,{\mathbb R})$ and
$\alpha\in {\mathbb R}^{(1, 2)}.$

\medskip

We know that $SL(2,{\mathbb R})$ acts on ${{\mathbb H}}$
transitively by
$$g<\tau>=(a\tau+b)(c\tau+d)^{-1},\ \ \
\begin{pmatrix} a & b \\ c & d \end{pmatrix}\in SL(2,{\mathbb R}),\ \
\tau\in{\mathbb H}.$$
Now we observe that the action (1.2) of $G$
on ${\mathbb H}\times{\mathbb C}$ may be rewritten as
$$(g,\alpha)\circ
(\tau,z)=\,(\,^t\!g^{-1}<\tau>, \,(z+\alpha_1\tau+\alpha_2)
(-b\tau+a)^{-1}\,),$$ where $g=\begin{pmatrix}
a&b\\c&d\end{pmatrix}\in SL(2,{\mathbb R}),\
\alpha=(\alpha_1,\alpha_2)\in {\mathbb R}^{(1, 2)},$ and
$(\tau,z)\in{\mathbb H}\times{\mathbb C}.$ Since the action (1.2)
is transitive and $K$ is the stabilizer of this action at the
origin $(i,0),\ {\mathbb H}\times{\mathbb C}$ can be identified
with the homogeneous space $G/K$ as follows\,:
\begin{equation}
G/K\ni (g,\alpha)K\longmapsto (g,\alpha)\circ (i,0).
\end{equation}
We see that we can express an element $Y$ of $S{{\mathcal P}}_2$
uniquely as
\begin{equation}
Y=\begin{pmatrix} y^{-1} & 0\\ 0& y
\end{pmatrix} \left[\begin{pmatrix} 1 & -x\\ 0 & 1\end{pmatrix}\right]=\,
\begin{pmatrix} y^{-1} & -xy^{-1}\\ -xy^{-1} & x^2y^{-1}+y\end{pmatrix}
\end{equation}
with $x,y\in {\mathbb R}$ and $y>0.$
\vspace{0.1in}\\
\noindent {\bf Lemma 2.1.} {\it We define the mapping $T:
S{\mathcal P}_2\times {\mathbb R}^{(1,2)}\longrightarrow {\mathbb
H}\times{\mathbb C}$ by
\begin{equation}
T(Y, V)=\,(x+iy,\,v_1(x+iy)+v_2\,), \end{equation} where $Y$ is of
the form {\rm (2.4)} and $V=(v_1,v_2)\in {\mathbb R}^{(1, 2)}.$
Then the mapping $T$ is a bijection which is compatible with the
above two actions {\rm (1.2)} and} (2.1).

\medskip

For any $Y\in S{{\mathcal P}}_2$ of the form (2.4), we put
\begin{equation}
g_Y=\,\begin{pmatrix} 1 & 0 \\ -x & 1\end{pmatrix}
\begin{pmatrix} y^{-1/2} & 0 \\ 0 & y^{1/2}\end{pmatrix}=\,
\begin{pmatrix} y^{-1/2} & 0\\ -xy^{-1/2} & y^{1/2}\end{pmatrix}.
\end{equation}
and \begin{equation} \alpha_{Y,V}=\,V\,\,^tg^{-1}_Y.
\end{equation}
Then we have
\begin{equation}
T(Y,V)=\,(\,g_Y,\,\alpha_{Y,V}\,)\circ
(i,0).
\end{equation}
\vspace{0.05in}\\
\noindent {\it Proof.} It is easy to prove the lemma. So we leave
the proof to the reader. \hfill$\square$

\bigskip

\indent Now we give a complete description of the algebra
$\mathbb{D}(\mathbb{H} \times {\mathbb C})$ of all differential
operators on ${\mathbb H}\times{\mathbb C}$ invariant under the
action (1.2) of $G$. First we note that the Lie algebra
${\mathfrak g}$ of $G$ is given by ${\mathfrak
g}=\,\left\{\,(X,Z)\,\vert\ X\in {\mathbb R}^{(2,2)},\
\sigma(X)=0,\ Z\in {\mathbb R}^{(1,2)} \,\right\}$ equipped with
the following Lie bracket
$$[(X_1,Z_1),\,(X_2,Z_2)]=\,(\,[X_1,X_2]_0,\,Z_2\,^tX_1-Z_1\,^t\!X_2\,),$$
where $[X_1,X_2]_0=\,X_1X_2-X_2X_1$ denotes the usual matrix
bracket and $(X_1,Z_1),\,(X_2,$

\par

\noindent $Z_2)\in {\mathfrak g}.$ And ${\mathfrak g}$ has the
following decomposition
$${\mathfrak g}={\mathfrak k}\oplus{\mathfrak p}\ \ \ \ \ \ \text{(\,direct\ sum\,)},$$
where ${\mathfrak k}=\,\bigg\{\,(X,0)\in {\mathfrak g}\ \vert\
X=\begin{pmatrix} 0 & x \\ -x & 0\end{pmatrix},\ \ x\in {\mathbb
R}\,\bigg\}$ and ${\mathfrak p}=\,\bigg\{\,(X,Z)\in{\mathfrak g}
\vert\ X=\,^tX\in {\mathbb R}^{(2,2)},\ \sigma(X)=0,\ Z\in
{\mathbb R}^{(1,2)}\,\bigg\}.$ \\
We observe that ${\mathfrak k}$ is the
Lie algebra of $K$ and that we have the following relations
$$[{\mathfrak k},\,{\mathfrak k}]\subset {\mathfrak k} \ \ \ \ \ \ \ \text{and}\ \ \ \ \ \ \
[{\mathfrak k},\,{\mathfrak p}]\subset {\mathfrak p}.$$ Thus the
coset space $G/K\cong {\mathbb H}\times{\mathbb C}$ is a {\it
reductive} homogeneous space in the sense of [12], \,p.\,284. It
is easy to see that the adjoint action $\text{Ad}$ of $K$ on
${\mathfrak p}$ is given by
\begin{equation}
\text{Ad}\,(k)((X,Z))\,=\,(\,kX\,^tk,\,Z\,^tk\,),
\end{equation}
where $k\in K$ and $(X,Z)\in {\mathfrak p}$ with $X=\,^t\!X,\
\sigma(X)=0.$ The action (2.9) extends uniquely to the action
$\rho$ of $K$ on the polynomial algebra $\text{Pol}\,({\mathfrak
p})$ of ${\mathfrak p}$ given by
\begin{equation}
\rho\,:\,K\longrightarrow \text{Aut\,(\,Pol}({\mathfrak p}\,)).
\end{equation}
Let $\text{Pol}\,({\mathfrak p})^K$ be the subalgebra of
$\text{Pol}\,({\mathfrak p})$ consisting of all invariants of the
action $\rho$ of $K$. Then according to [12], \,Theorem 4.9,\,
p.\,287, there exists a canonical linear bijection
$\lambda\,(\,P\longmapsto D_{\lambda(P)})$ of
$\text{Pol}\,({\mathfrak p})^K$ onto $\mathbb{D}(\mathbb{H} \times
{\mathbb C})$. Indeed, if $(\xi_k)\,(\,1\leq k\leq 4\,)$ is any
basis of ${\mathfrak p}$ and $P\in \text{Pol}\,({\mathfrak p})^K,$
then
\begin{equation}
\left(\,D_{\lambda(P)}f\right)({\tilde g}\circ (i,0))\,=\,
\left[\,P\left({{\partial}\over{\partial t_k}}\right)f( ({\tilde
g} \ast \text{exp}\,(\sum_{k=1}^4t_k\xi_k))\circ (i,0))
\right]_{(t_k)=0},\end{equation} where ${\tilde g}\in G$ and $f\in
C^{\infty}({\mathbb H}\times{\mathbb C}).$

\medskip

We put $$e_1=\,\left( \begin{pmatrix} 1 & 0\\ 0 &
-1\end{pmatrix},(0,0)\right),\ \ \ e_2=\,\left( \begin{pmatrix} 0
& 1 \\ 1
& 0\end{pmatrix},(0,0)\right)$$ and $$ f_1=\,\left( \begin{pmatrix} 0 & 0\\
0 & 0\end{pmatrix},(1,0)\right),\ \ \ f_2=\,\left( \begin{pmatrix} 0 & 0\\
0 & 0\end{pmatrix},(0,1)\right).$$ Then $e_1,\,e_2,\,f_1,\,f_2$
form a basis of ${\mathfrak p}.$ We write for coordinates $(X,Z)$
by
$$X\,=\,\begin{pmatrix} x & y\\ y & -x\end{pmatrix}\ \ \ \text{and}\ \ \
Z\,=\,(z_1,z_2)$$ with real variables $x, y, z_1$ and $z_2.$
\vspace{0.1in}\\
\noindent {\bf Lemma 2.2.} {\it The following polynomials
\begin{eqnarray*}
P(X,Z) &=&{\frac 18}\,\sigma(X^2)\,=\,{\frac 14}\,(x^2+y^2),\\
\xi(X,Z) &= & Z\,^tZ\,=\,z_1^2+z_2^2,\\
P_1(X,Z) &=& -\,{\frac 12}\,ZX\,^tZ\,=\,{\frac
12}\,(z_2^2-z_1^2)x\,-\, z_1z_2 y ~~~{\text
and} \\
P_2(X,Z) &=& {\frac 12}\,(z_2^2-z_1^2)\,y\,+\,z_1z_2 \,x
\end{eqnarray*}
are algebraically independent generators of
$\text{Pol}\,({\mathfrak p})^K.$}
\vspace{0.05in}\\
\noindent {\it Proof.} We leave the proof of the above lemma to
the reader. \hfill$\square$

\medskip

Now we are ready to compute the $G$-invariant differential
operators $D,\,\Psi,\,D_1$ and $D_2$ corresponding to the
$K$-invariants $P,\,\xi,\,P_1$ and $P_2$ respectively under the
canonical linear bijection (2.11). For real variables
$t=(t_1,t_2)$ and $s=(s_1,s_2),$ we have
$$\exp \,(\,t_1e_1+t_2e_2+s_1f_1+s_2f_2\,) \,\ =\,\ \Bigg(\,\begin{pmatrix} a_1(t,s) & a_3(t,s)\\ a_3(t,s)
a_2(t,s) \end{pmatrix},\,\,(\,b_1(t,s),\,b_2(t,s)\,\Bigg),$$ where
\begin{align*}
 a_1(t,s)&=\,1\,+\,t_1\,+\,{1\over
{2!}}\,(\,t_1^2+t_2^2\,)\,+\, {1\over{3!}}\,t_1(t_1^2+t_2^2\,)
\,+\,{1\over {4!}}\,(t_1^2+t_2^2)^2\,+\,\cdots\,\\
a_2(t,s)&=\,1\,-\,t_1\,+\,{1\over {2!}}\,(\,t_1^2+t_2^2\,)\,-\,
{1\over {3!}}\,t_1(t_1^2+t_2^2\,)
\,+\,{1\over {4!}}\,(t_1^2+t_2^2)^2\,-\,
\cdots\,,\\
a_3(t,s)&=\,t_2\,+\,{1\over {3!}}\,t_2(\,t_1^2\,+\,t_2^2\,)\,+\,
{1\over {5!}}\,t_2(\,t_1^2\,+\,t_2^2\,)^2\,+\,\cdots\,,\\
b_1(t,s)&=\,s_1\,-\,{1\over {2!}}\,(s_1t_1+s_2t_2)\,+\,
{1\over {3!}}\,s_1(t_1^2+t_2^2)\,-\,{1\over {4!}}\,(s_1t_1+s_2t_2)
(t_1^2+t_2^2)\,+\cdots\,,\\
b_2(t,s)&=\,s_2\,-\,{1\over {2!}}\,(s_1t_2-s_2t_1)\,+\,{1\over{3!}}\,
s_2(t_1^2+t_2^2)\,-\,{1\over {4!}}\,(s_1t_2-s_2t_1)(t_1^2+t_2^2)\,
+\,\cdots\,.
\end{align*}
For brevity, we write $a_j,\,b_k$ for $a_j(t,s),\,b_k(t,s)\
(\,j=1,2,3,\,k=1,2\,)$ respectively. We now fix an element
$(g,\alpha)\in G$ and write
$$g\,=\,\begin{pmatrix} g_1 & g_{12}\\ g_{21} & g_2 \end{pmatrix} \in SL(2,{\mathbb R})
\ \ \ \text{and}\ \ \ \alpha\,=\,(\alpha_1,\alpha_2)\in {\mathbb
R}^{(1,2)}.$$ We put $\left(\,\tau(t,s),\,z(t,s)\,\right)=
\,\left(\,(g,\alpha)\ast\text{exp}\,
\left(\,t_1e_1\,+\,t_2e_2\,+\,s_1f_1\,+\,s_2f_2\,\right)\right)\circ
(i,0)$ with $\tau(t,s)=\,x(t,s)\,+\,i\,y(t,s)\ \ \ \text{and}\ \ \
z(t,s)=\,u(t,s)\,+\,i\,v(t,s).$\\
Here $x(t,s),\,y(t,s),\,u(t,s)$ and $v(t,s)$ are real. By an easy
calculation, we obtain
\begin{align*}
x(t,s)\,&=\,-(\,\tilde {a}\tilde {c}+\tilde {b}\tilde {d}\,)\,(\,\tilde {a}^2+\tilde {b}^2\,)^{-1},\\
y(t,s)\,&=\,(\,\tilde {a}^2+\tilde {b}^2\,)^{-1},\\
u(t,s)\,&=\,(\,\tilde {a}\,{\tilde {\alpha}}_2-\tilde {b}{\tilde {\alpha}}_1\,)\,(\,\tilde {a}^2+\tilde {b}^2\,)^{-1},\\
v(t,s)\,&=\,(\,\tilde {a}\,{\tilde {\alpha}}_1+\tilde {b}\,{\tilde
{\alpha}}_2\,)\,(\,\tilde {a}^2+\tilde {b}^2\,)^{-1},
\end{align*}
where $\tilde {a} =\,g_1a_1\,+\,g_{12}a_3,$ \,\ $\tilde {b}
=\,g_1a_3\,+\,g_{12}a_2,$\,\ $\tilde {c} =\,g_{21}a_1\,+\,g_2a_3,$
\,\ $\tilde {d} =\,g_{21}a_3\,+\,g_2a_2,$\,\ ${\tilde {\alpha}}_1
=\,\alpha_1a_2\,-\,\alpha_2a_3\,+\,b_1,$\,\ ${\tilde {\alpha}}_2
=\,-\alpha_1a_3\,+\,\alpha_2a_1\,+\,b_2.$ \\
By an easy calculation, at $t=s=0$, we have
\begin{align*}
{{\partial x}\over {\partial t_1}}\,&=\,4\,g_1\,g_{12}\,
(\,g_1^2\,+\,g_{12}^2\,)^{-2},\\
{{\partial y}\over {\partial t_1}}\,&=\,-2\,(\,g_1^2\,-\,g_{12}^2\,)\,
(\,g_1^2\,+\,g_{12}^2\,)^{-2},\\
{{\partial u}\over {\partial t_1}}\,&=\,4\,g_1\,g_{12}\,(\,g_1\,
\alpha_1\,+\,g_{12}\,\alpha_2\,)\,(\,g_1^2\,+\,g_{12}^2\,)^{-2},\\
{{\partial v}\over {\partial
t_1}}\,&=\,-\,2\,(\,g_1\,\alpha_1\,+\,
g_{12}\,\alpha_2\,)\,(\,g_1^2\,-\,g_{12}^2\,)\,(\,g_1^2\,+
\,g_{12}^2\,)^2,\\
{{\partial^2 x}\over {\partial t_1^2}}\,&=\,-\,16\,g_1\,g_{12}\,
(\,g_1^2\,-\,g_{12}^2\,)\,(\,g_1^2\,+\,g_{12}^2\,)^{-3},\\
{{\partial^2 y}\over {\partial t_1^2}}\,&=\,8\,(\,g_1^2\,-\,g_{12}^2
\,)^2\,(\,g_1^2\,+\,g_{12}^2\,)^{-3}\,-\,4\,(\,g_1^2\,+\,
g_{12}^2\,)^{-1},\\
{{\partial^2 u}\over {\partial t_1^2}}\,&=\,-\,16\,g_1\,g_{12}\,
(\,g_1\,\alpha_1\,+\,g_{12}\,\alpha_2\,)\,(\,g_1^2\,-\,g_{12}^2\,)\,
(\,g_1^2\,+\,g_{12}^2\,)^{-3},\\
{{\partial^2 v}\over {\partial t_1^2}}\,&=\,4\,
(\,g_1\,\alpha_1\,+\,g_{12}\,\alpha_2\,)\,(\,g_1^4\,+\,g_{12}^4\,-\,
6\,g_1^2\,g_{12}^2\,)\,(\,g_1^2\,+\,g_{12}^2\,)^{-3}
\end{align*}
and
\begin{align*}
{{\partial x}\over {\partial
t_2}}\,&=\,-\,2\,
(\,g_1^2\,-\,g_{12}^2\,)\,(\,g_1^2\,+\,g_{12}^2\,)^{-2},\\
{{\partial y}\over {\partial t_2}}\,&=\,-\,4\,g_1\,g_{12}\,
(\,g_1^2\,+\,g_{12}^2\,)^{-2},\\
{{\partial u}\over {\partial t_2}}\,&=\,-\,2\,
(\,g_1\,\alpha_1\,+\,g_{12}\,\alpha_2\,)\,(\,g_1^2\,-\,g_{12}^2\,)\,
(\,g_1^2\,+\,g_{12}^2\,)^{-2},\qquad\qquad
\end{align*}
\begin{align*}
{{\partial v}\over {\partial t_2}}\,&=\,-\,4\,g_1\,g_{12}\,
(\,g_1\,\alpha_1\,+\,g_{12}\,\alpha_2\,)\,(\,g_1^2\,+\,
g_{12}^2\,)^{-2},\\
{{\partial^2 x}\over {\partial t_2^2}}\,&=\,16\,g_1\,g_{12}\,
(\,g_1^2\,-\,g_{12}^2\,)\,(\,g_1^2\,+\,g_{12}^2\,)^{-3},\\
{{\partial^2 y}\over {\partial t_2^2}}\,&=\,32\,
g_1^2\,g_{12}^2\,(\,g_1^2\,+\,g_{12}^2\,)^{-3}\,-\,
4\,(\,g_1^2\,+\,g_{12}^2\,)^{-1},\\ {{\partial^2 u}\over {\partial
t_2^2}}\,&=\,16\,g_1\,g_{12}\,
(\,g_1\,\alpha_1\,+\,g_{12}\,\alpha_2\,)\,(\,g_1^2\,-\,g_{12}^2\,)\,
(\,g_1^2\,+\,g_{12}^2\,)^{-3},\\ {{\partial^2 v}\over {\partial
t_2^2}}\,&=\,-\,4\,
(\,g_1\,\alpha_1\,+\,g_{12}\,\alpha_2\,)\,(\,g_1^4\,+\,g_2^4\,-\,
6\,g_1\,g_{12}^2\,)\,(\,g_1^2\,+\,g_{12}^2\,)^{-3}.
\end{align*}
We note that $\tilde {a}\tilde {d}\,-\,\tilde {b}\tilde
{c}\,=\,1,\ \ a_1a_2-a_3^2=1\ \ \text{and}\ \
g_1g_2-g_{12}g_{21}=1.$\\
\indent Using the above facts and applying the chain rule, we can
easily compute the differential operators $D,\,\Psi,\,D_1$ and
$D_2.$ It is known that the images of generators $P,\,\xi,\,P_1$
and $P_2$ under $\lambda$ are generators of $\mathbb{D}(\mathbb{H}
\times \mathbb{C})$\,(\,cf.\,[11]).

\medskip

Summarizing, we have the following.
\vspace{0.1in}\\
\noindent {\bf Theorem 2.3.}  {\it The algebra
$\mathbb{D}(\mathbb{H} \times \mathbb{C})$ is generated by the
following differential operators
\begin{align}
D=\,y^2\,\left(\,{{\partial^2}\over {\partial x^2}}\,+\,
{{\partial^2}\over {\partial y^2}}\,\right)
\,+\,v^2\,\left(\,{{\partial^2}\over{\partial
u^2}}\,+\,{{\partial^2}\over{\partial v^2}}\,\right)
+\,2\,y\,v\,\left(\,{{\partial^2}\over{\partial x\partial
u}}\,+\,{{\partial^2}\over{\partial y\partial v}}\,\right),
\end{align}
\begin{equation}
\Psi=\,y\,\left(\,{{\partial^2}\over {\partial u^2}}\,+\,
{{\partial^2}\over {\partial v^2}}\,\right),
\end{equation}
\begin{equation}
D_1=\,2\,y^2\,{{\partial^3}\over {\partial x\partial u
\partial v}}\,-\,y^2\,{{\partial}\over{\partial y}}
\left(\,{{\partial^2}\over{\partial u^2}}\,-\,
{{\partial^2}\over{\partial v^2}}\,\right)\,+\,\left(\,
v\,{{\partial}\over{\partial v}}\,+\,1\,\right)\Psi
\end{equation}
and
\begin{equation}
D_2=\,y^2\,{{\partial}\over{\partial x}}\left(\,
{{\partial^2}\over{\partial v^2}}\,-\,{{\partial^2}\over {\partial
u^2}}\,\right)\,-\,2\,y^2\,{{\partial^3}\over{\partial y\partial u
\partial v}}\,-\,v\,{{\partial}\over{\partial u}}\Psi,
\end{equation}
where $\tau=x+iy$ and $z=u+iv$ with real variables $x,y,u,v.$
Moreover, we have
\begin{align*}
[D,\,\Psi]=D\Psi-\Psi D\,=&\,
2\,y^2\,{{\partial}\over{\partial
y}}\left(\,{{\partial^2}\over{\partial
u^2}}\,-\,{{\partial^2}\over{\partial v^2}}\,\right)\,-\,
4\,y^2\,{{\partial^3}\over{\partial x\partial u\partial v}}\\
&\ \ \,-\,2\,\left(\,v\,{{\partial}\over{\partial
v}}\Psi\,+\,\Psi\,\right).
\end{align*}
In particular, the algebra $\mathbb{D}(\mathbb{H} \times {\mathbb
C})$ is not commutative. Thus the homogeneous space ${{\mathbb
H}}\times {\mathbb C}$ is not weakly symmetric in the sense of A.
Selberg} ([19]).

\bigskip

Now we provide a natural $G$-invariant K{\"a}hler metric on
${\mathbb H}\times{\mathbb C}.$
\vspace{0.1in}\\
\noindent {\bf Proposition 2.4.}  {\it The Riemannian metric
$ds^2$ on ${\mathbb H}\times{\mathbb C}$ defined by
$$ ds^2\,=\,{{y\,+\,v^2}\over {y^3}}\,(\,dx^2\,+\,dy^2\,)\,+\, {\frac
1y}\,(\,du^2\,+\,dv^2\,) \ -\,{{2v}\over {y^2}}\,
(\,dx\,du\,+\,dy\,dv\,)$$
is invariant under the action {\rm
(1.2)} of $G$ and is a K{\"a}hler metric on ${{\mathbb H}}\times
{\mathbb C}$. The Laplace-Beltrami operator $\Delta$ of the
Riemannian space $(\,{\mathbb H}\times{\mathbb C},\,ds^2\,)$ is
given by
$$\Delta\,=\, y^2\,\left(\,{{\partial^2}\over{\partial
x^2}}\,+\,{{\partial^2}\over{\partial y^2}}\,\right)\,+\,
(\,y\,+\,v^2\,)\,\left(\,{{\partial^2}\over{\partial
u^2}}\,+\,{{\partial^2}\over{\partial v^2}}\,\right) \
+\,2\,y\,v\,\left(\,{{\partial^2}\over{\partial x\partial
u}}\,+\,{{\partial^2}\over{\partial y\partial v}}\,\right).
$$
That is, $\Delta\,=\,D\,+\,\Psi.$}
\vspace{0.05in}\\
\noindent {\it Proof.} For $Y\in S{{\mathcal P}}_2$ of the form
(2.4) and $(v_1, v_2)\in {\mathbb R}^{(1,2)},$ it is easy to see
that
$$dY=\,\begin{pmatrix}
-\,y^{-2} \,dy & -\,y^{-1}\,dx\,+\,x\,y^{-2}\,dy\\
-\,y^{-1}\,dx\,+\,x\,y^{-2}\,dy & 2\,
x\,y^{-1}\,dx\,+\,(\,1\,-\,x^2\,y^{-2}\,)\,dy
\end{pmatrix}$$
and $dV=(dv_1, dv_2).$ Then we can show that the following metric
$d{\tilde s}^2$ on $S{\mathcal P}_2 \times {\mathbb R}^{(1,2)}$
defined by
$$d{\tilde s}^2\,=\,{{dx^2\,+\,dy^2}\over {y^2}}\,+\,{\frac 1y}\,
\left\{\,(\,x^2\,+\,y^2\,)\,dv_1^2\,+\,2\,x\,dv_1\,dv_2\,+\,
dv_2^2\,\right\}$$
is invariant under the action (2.1) of $G$.
Indeed, since
$$Y^{-1}\,=\,\begin{pmatrix} y\,+\,x^2\,y^{-1} &
x\,y^{-1}\\ x\,y^{-1} & y^{-1}\end{pmatrix},$$ we can easily show
that $d{\tilde s}^2\,=\,{\frac 12}\,\sigma(Y^{-1}dY
Y^{-1}dY)\,+\,dV\,Y^{-1}\,^t(dV).$\\
For an element $(\,g,\,\alpha\,)\in G$ with $g\in SL(2,{\mathbb
R})$ and $\alpha\in {\mathbb R}^{(1,2)},$ we put
$$(\,Y^{\ast},\,V^{\ast}\,)=\,(\,g,\,\alpha\,)\cdot (Y,V)\,=\,
(\,gY\,^tg,\,(V+\alpha)\,^tg\,).$$ Since $Y^{\ast}\,=\,gY\,^tg \,\
\text{and} \,\ V^{\ast}\,=\,(\,V\,+\,\alpha\,) \,^tg,$ we get
$dY^{\ast}\,=\,g\,dY\,^tg \,\ \text{and}\,\
V^{\ast}\,=\,(\,V\,+\,\alpha\,)\,^tg.$\\
Therefore by a simple calculation, we can show that
\begin{align*}
&\
\sigma\left(\,Y^{\ast -1}\,dY^{\ast}\,Y^{\ast
-1}\,dY^{\ast}\,\right) \,+\,dV^{\ast}\,Y^{\ast
-1}\,^t(dV^{\ast})\\ &=\,
\sigma(\,Y^{-1}dY\,Y^{-1}\,dY\,)\,+\,dV\,Y^{-1}\,^t(dV).
\end{align*}
Hence the metric $d{\tilde s}^2$ is invariant under the action
(2.1) of $G$.

\medskip

Using this fact and Lemma 2.1, we can prove that the metric $ds^2$
in the above theorem is invariant under the action (1.2). Since
the matrix form $(\,g_{ij}\,)$ of the metric $ds^2$ is given by
$$(\,g_{ij}\,)\,=\,\begin{pmatrix} (\,y\,+\,v^2\,)\,y^{-3} & 0 &
-\,v\,y^{-2} & 0 \\ 0 & (\,y\,+\,v^2\,)\,y^{-3} & 0 &
-\,v\,y^{-2}\\ -\,v\,y^{-2} & 0 & y^{-1} & 0 \\ 0 & -\,v\,y^{-2} &
0 & y^{-1} \end{pmatrix}$$ and $\det\,(\,g_{ij}\,)\,=\,y^{-6},$
the inverse matrix $(\,g^{ij}\,)$ of $(\,g_{ij}\,)$ is easily
obtained by
$$(\,g^{ij}\,)\,=\,\begin{pmatrix} y^2 & 0 & y\,v & 0 \\ 0 & y^2 & 0 &
y\,v \\ y\,v & 0 & y\,+\,v^2 & 0 \\ 0 & y\,v & 0 & y\,+\,v^2
\end{pmatrix}.$$ Now it is easily shown that $D\,+\,\Psi$ is the
Laplace-Beltrami operator of $(\,{\mathbb H}\times{\mathbb
C},\,ds^2\,).$ \hfill $\square$
\vspace{0.1in}\\
\noindent {\bf Remark 2.5.} We can show that for any two positive
real numbers $\alpha$ and $\beta$, the following metric
$$ds_{\alpha,\beta}^2=\,\alpha\,{{dx^2\,+\,dy^2}\over {y^2} }\,+\,\beta\, {
{v^2(dx^2\,+\,dy^2)\,+\,y^2(du^2\,+\,dv^2)\,-\,2\,yv\,(dx\,du\,+\,
dy\,dv)} \over {y^3} }$$ is also a Riemannian metric on ${\mathbb
H}\times{\mathbb C}$ which is invariant under the action (1.2) of
$G$. In fact, we can see that the two-parameter family of
$ds^2_{\alpha,\beta}\,(\,\alpha>0,\ \beta>0\,)$ provides a
complete family of Riemannian metrics on ${\mathbb
H}\times{\mathbb C}$ invariant under the action of (1.2) of $G.$
It can be easily seen that the Laplace-Beltrami operator
$\Delta_{\alpha,\beta}$ of $ds^2_{\alpha,\beta}$ is given by
\begin{align*}
 \Delta_{\alpha,\beta}\,&=\,{ {1}\over {\alpha}
}\,y^2\,(\,{{\partial^2}\over{\partial
x^2}}\,+\,{{\partial^2}\over{\partial y^2}}\,)\,+\, \left(\,{
{y}\over{\beta} }\,+\,{
{v^2}\over {\alpha} }\,\right)\,\left(\, {{\partial^2}\over{\partial u^2}}\,+\,{{\partial^2}\over{\partial v^2}}\,\right)\\
&\quad\ \ +\, { {2\,yv}\over {\alpha} }\,\left(\,
{{\partial^2}\over{\partial x\partial
u}}\,+\,{{\partial^2}\over{\partial y\partial v}}\,\right)\\ &=\ {
{1}\over {\alpha} }\,D\,+\,{ {1}\over {\beta} }\,\Psi.
\end{align*}
\vspace{0.1in}\\
\noindent {\bf Remark 2.6.}  By a tedious computation, we see that
the scalar curvature of $({{\mathbb H}}\times {\mathbb C},\,ds^2)$
is $-3.$

\medskip

We want to propose the following problem to be studied in the
future.
\vspace{0.1in}\\
\noindent {\bf Problem 2.7.} Find all the eigenfunctions of
$\Delta$.

\medskip

We will give some examples of eigenfunctions of $\Delta$.
\begin{itemize}
\item[(1)] $h(x,y)=y^{1\over 2}K_{s-{\frac12}}(2\pi |a|y)\,e^{2\pi
iax}\quad (s\in {\mathbb C},\ a\not=0\,)$ with eigenvalue
$s(s-1),$

\smallskip

\noindent where
\begin{equation}
K_s(z):={\frac12}\int^{\infty}_0 \exp\left\{-{z\over
2}(t+t^{-1})\right\}\,t^{s-1}\,dt,\quad{\text Re} \, z >
0.\end{equation}
\item[(2)] $y^s,\ y^s x,\ y^s u\ (s\in{\mathbb
C})$ with eigenvalue $s(s-1).$
\item[(3)] $y^s v,\ y^s uv,\ y^s xv$ with eigenvalue $s(s+1).$
\item[(4)] $x,\,y,\,u,\,v,\,xv,\,uv$ with eigenvalue $0$.
\item[(5)] All Maass wave forms.
\end{itemize}
\vspace{0.2in}
\noindent{\bf 3. Maass-Jacobi forms}
\setcounter{equation}{0}
\renewcommand{\theequation}{3.\arabic{equation}}
\vspace{0.1in}\\
\indent Let $\Delta$ be the Laplace-Beltrami operator of the
Riemannian metric $ds^2$ on ${\mathbb H}\times{\mathbb C}$ defined
in Proposition 2.4. Using this operator, we define the notion of
Maass-Jacobi forms.
\vspace{0.1in}\\
\noindent {\bf Definition 3.1.} A smooth bounded function
$f:{\mathbb H}\times{\mathbb C} \longrightarrow {\mathbb C}$ is
called a {\it Maass}-{\it Jacobi\ form} if it satisfies the
following conditions (MJ1)-(MJ3)\,:
\begin{itemize}
\item[(MJ1)] $f({\tilde {\gamma}}\circ
(\tau,z))\,=\,f(\tau,z) \ \ \ \text{for\ all}\ {\tilde
{\gamma}}\in \Gamma$ and $(\tau,z)\in {\mathbb H}\times{\mathbb
C}$.
\item[(MJ2)] $f$ is an eigenfunction of the Laplace-Beltrami
operator $\Delta.$
\item[(MJ3)] $f$ has a polynomial growth, that
is, $f$ fulfills a boundedness condition.
\end{itemize}

For a complex number $\lambda\in {\mathbb C}$, we denote by
$MJ(\Gamma,\lambda)$ the vector space of all Maass-Jacobi forms
$f$ such that $\Delta f=\lambda f.$ We note that, since $\Delta
f=\lambda f$ is an elliptic partial differential equation,
Maass-Jacobi forms are real analytic (see [8]). Professor Berndt
kindly informed me that he also considered such automorphic forms
in ([1]) (also see [4], p.82).

\medskip

Let $f\in MJ(\Gamma,\lambda)$ be a Maass-Jacobi form with
eigenvalue $\lambda.$ Then it is easy to see that the function
$\phi_f:G\longrightarrow {\mathbb C}$ defined by
\begin{equation}
\phi_f(g,\alpha)=f((g,\alpha)\circ (i,0)),\ \ \ (g,\alpha)\in G
\end{equation}
satisfies the following conditions (MJ$1)^0$-(MJ$3)^0$:
\begin{itemize}
\item[(MJ$1)^0$] $\phi_f(\gamma xk)=\phi_f(x)\ \ \text{for\
all}\ \gamma\in \Gamma,\ x\in G\ \text{and}\ k\in K.$
\item[(MJ$2)^0$] $\phi_f$ is an eigenfunction of the
Laplace-Beltrami operator $\Delta_0$ of $(G, \,ds_0^2),$ where
$ds_0^2$ is a $G$-invariant Riemannian metric on $G$ induced by
$({\mathbb H}\times{\mathbb C},\,ds^2).$
\item[(MJ$3)^0$] $\phi_f$ has a suitable polynomial growth\,(cf.\,[5]).
\end{itemize}

For any right $K$-invariant function $\phi:G\longrightarrow
{\mathbb C}$ on $G$, we define the function $f_{\phi}:{\mathbb
H}\times{\mathbb C}\longrightarrow {\mathbb C}$ by
\begin{equation}
f_{\phi}(\tau,z)=\phi(g,\alpha),\ \ \ (\tau,z)\in {\mathbb
H}\times{\mathbb C},
\end{equation}
where $(g,\alpha)$ is an element of $G$ such that $(g,
\alpha)\circ (i,0)=(\tau, z).$ Obviously it is well defined
because (3.2) is independent of the choice of $(g,\alpha)\in G$
such that $(g, \alpha)\circ (i, 0)=(\tau,z).$ It is easy to see
that if $\phi$ is a smooth bounded function on $G$ satisfying the
conditions (MJ$1)^0$-(MJ$3)^0,$ then the function $f_{\phi}$
defined by (3.2) is a Maass-Jacobi form.

\medskip

Now we characterize Maass-Jacobi forms as smooth eigenfunctions on
$S{\mathcal P}_n\times {\mathbb R}^{(m, n)}$ satisfying a certain
invariance property.
\vspace{0.1in}\\
\noindent {\bf Proposition 3.2.} {\it Let $f:{\mathbb
H}\times{\mathbb C}\longrightarrow{\mathbb C}$ be a nonzero
Maass-Jacobi form in $MJ(\Gamma, \lambda).$ Then the function $h_f
: S{\mathcal P}_2\times {\mathbb R}^{
(1,2)}\longrightarrow{\mathbb C}$ defined by
\begin{equation}
h_f(Y,V)=f((g,\,V\,^tg^{-1})\circ (i,0))\ \ \text{for\ some}\ g\in
SL(2,{\mathbb R})\ \text{with}\ Y=g\,^tg \end{equation} satisfies
the following conditions\,:
\begin{itemize}
\item [{\rm (MJ$1)^{\ast}$}] $h_f(\gamma
Y\,^t\gamma,\,(V+\delta)\,^t\gamma)=h_f(Y,V)\ \ \ \ \text{for\
all}\ (\gamma,\delta)\in \Gamma\ \text{with}\ \gamma\in
SL(2,{\mathbb Z}) \ \text{and}\ \delta\in {\mathbb Z}^{(1,2)}.$
\item [{\rm (MJ$2)^{\ast}$}] $h_f$ is an eigenfunction of the
Laplace-Beltrami operator ${\tilde \Delta}$ on the homogeneous
space $(S{\mathcal P}_2\times {\mathbb R}^{(1,2)},d{\tilde s}^2)$,
where $d{\tilde s}^2$ is the $G$-invariant Riemannian metric on
$S{\mathcal P}_2\times {\mathbb R}^{(1,2)}$ induced from $d{\tilde
s}^2.$
\item[{\rm (MJ$3)^{\ast}$}] $h_f$ has a suitable
polynomial growth.
\end{itemize}

Here if $(Y, V)$ is a coordinate of $S{\mathcal P}_2\times
{\mathbb R}^{(1,2)}$ given in Lemma 2.1, then the $G$-invariant
Riemannian metric $d{\tilde s}^2$ and its Laplace-Beltrami
operator ${\tilde \Delta}$ on $S{\mathcal P}_2\times {\mathbb
R}^{(1,2)}$ are given by
$$d{\tilde s}^2={1\over {y^2}}(dx^2+dy^2)\,+\,{1\over y}\left\{
(x^2+y^2)dv_1^2\,+\,2xdv_1dv_2\,+\,dv_2^2\right\}$$
and
$${\tilde \Delta}=\,y^2\left({{\partial^2}\over{\partial x^2}}\,+\,
{{\partial^2}\over{\partial y^2}}\right)\,+\,{1\over y}\left\{
{{\partial^2}\over{\partial v_1^2}}-2x{{\partial^2}\over {\partial v_1
\partial v_2}}+(x^2+y^2){{\partial^2}\over{\partial v_2^2}}\right\}.$$
Conversely, if $h$ is a smooth bounded function on $S{\mathcal
P}_2\times {\mathbb R}^{(1,2)}$ satisfying the above conditions
{\rm (MJ$1)^{\ast}$-(MJ$3)^{\ast},$} then the function
$f_h:{\mathbb H}\times{\mathbb C}\longrightarrow{\mathbb C}$
defined by
\begin{equation}
f_h(\tau,z)=\,h(g\,^tg,\,\alpha\,^tg) \end{equation} for some $(g,
\alpha)\in G$ with $(g,\alpha)\circ (i, 0)=(\tau, z)$ is a
Maass-Jacobi form on ${\mathbb H}\times{\mathbb C}$.}
\vspace{0.05in}\\
\noindent {\it Proof.} First of all, we note that $h_f$ is well
defined because (3.3) is independent of the choice of $g$ with
$Y=\,g\,^tg.$ If $(\gamma,\delta)\in \Gamma$ with $\gamma\in
\Gamma_1,\ \delta\in {\mathbb Z}^{(1, 2)}$ and $(Y,V)\in
S{\mathcal P}_2\times {\mathbb R}^{(1,2)}$ with $Y=\,g\,^tg$ for
some $g\in SL(2, {\mathbb R}),$ then the element
$g_{\gamma}:=\gamma g$ satisfies $\gamma Y\,^t\gamma=\gamma
g\,^t(\gamma g).$\\
Thus according to the definition of $h_f$, for all
$(\gamma,\delta)\in \Gamma$ and $(Y,V)\in S{\mathcal P}_n\times
{\mathbb R}^{ (m, n)}$, we have
\begin{align*}
h_f(\gamma
Y\,^t\gamma,\,(V+\delta)\,^t\gamma)&=\,f((\gamma
g,\,(V+\delta)\,^t\gamma\, ^t(\gamma g)^{-1})\circ (i,0))\\
&=\,f((\gamma g,\,(V+\delta)\,^tg^{-1})\circ (i,0))\\
&=\,f(((\gamma,\delta)\ast (g,\,V\,^tg^{-1}))\circ (i,0))\\
&=\,f((g,\,V\,^tg^{-1})\circ (i,0))\ \ \ (\text{because}\ f\
\text{is}\ \Gamma\!-\!\text{invariant})\\ &=\,h_f(Y,V).
\end{align*}
Therefore this proves the condition (MJ$1)^{\ast}.\ d{\tilde s}^2$
and ${\tilde \Delta}$ are obtained from Lemma 2.1 and Proposition
2.3. Hence $h_f$ is an eigenfunction of ${\tilde \Delta}.$ Clearly
$h_f$ satisfies the condition (MJ$3)^{\ast}.$

\smallskip

Conversely we note that $f_h$ is well defined because (3.4) is
independent of the choice of $(g,\alpha)\in G$ with
$(g,\alpha)\circ (i,0)=(\tau,z).$ If ${\tilde
\gamma}=(\gamma,\delta)\in \Gamma$ and $(\tau,z)\in {\mathbb
H}\times{\mathbb C}$ with $(g,\alpha)\circ (i, 0)=(\tau,z ),$ then
we have
\begin{align*}
f_h({\tilde \gamma}\circ
(\tau,z))&=\,f_h({\tilde\gamma}\circ ((g,\alpha)\circ (i,0)))\\
&=\,f_h(({\tilde \gamma}\ast (g,\alpha))\circ
(i,0))\\ &=\,f_h((\gamma g,\,\delta\,^tg^{-1}+\alpha)\circ (i,0))\\
&=\,h((\gamma g)\,^t(\gamma g),\,(\delta\,^tg^{-1}+\alpha)\,^t(\gamma g))\\
&=\,h((\gamma (g\,^tg)\,^t\gamma,\,(\delta+\alpha\,^tg)\,^t\gamma)\\
&=\,h(g\,^tg,\,\alpha\,^tg)\\ &=\,f_h((g,\alpha)\circ
(i,0))=f_h(\tau,z).
\end{align*}
Thus $f_h$ satisfies the condition (MJ1). It is easy to see that
$f_h$ satisfies the conditions (MJ2) and (MJ3). \hfill $\square$
\vspace{0.1in}\\
\noindent {\bf Definition 3.3.} A smooth bounded function on $G$
or $S{\mathcal P}_2\times {\mathbb R}^{(1,2)}$ is also called a
{\it Maass}-{\it Jacobi\ form} if it satisfies the conditions
(MJ$1)^0$-(MJ$3)^0$ or (MJ$1)^{\ast}$-(MJ$3)^{\ast}$.
\vspace{0.1in}\\
\noindent {\bf Remark 3.4.} We note that Maass wave forms are
special ones of Maass-Jacobi forms. Thus the number of $\lambda$'s
with $MJ(\Gamma,\lambda)\neq 0$ is infinite.
\vspace{0.1in}\\
\noindent {\bf Theorem 3.5.} {\it For any complex number
$\lambda\in{\mathbb C},$ the vector space $MJ(\Gamma, \lambda)$ is
finite dimensional.}
\vspace{0.05in}\\
\noindent {\it Proof.} The proof follows from [10], Theorem 1,
p.\,8 and [5], \, p.\,191. \hfill $\square$
\vspace{0.2in}\\
\noindent{\bf 4. On the group $SL_{2, 1}({\mathbb R})$}
\setcounter{equation}{0}
\renewcommand{\theequation}{4.\arabic{equation}}
\vspace{0.1in}\\
\indent For brevity, we set $H={\mathbb R}^{(1, 2)}.$ Then we have
the split exact sequence
$$0\longrightarrow H\longrightarrow
G\longrightarrow SL(2,{\mathbb R})\longrightarrow 1.$$ We see that
the unitary dual ${\hat H}$ of $H$ is isomorphic to ${\mathbb
R}^2$. The unitary character $\chi_{(\lambda, \mu)}$ of $H$
corresponding to $(\lambda, \mu)\in {\mathbb R}^2$ is given by
$$\chi_{(\lambda,\mu)}(x,y)=e^{2\pi i(\lambda x+\mu y)},\quad\ \ (x,y)\in
H.$$ $G$ acts on $H$ by conjugation and hence this action induces
the action of $G$ on ${\hat H}$ as follows.
\begin{equation}
G\,\times\,{\hat H}\,\longrightarrow\,{\hat H},\quad\
(g,\chi)\mapsto \chi^g,\quad g\in G,\ \chi\in {\hat H},
\end{equation}
where the character $\chi^g$ is defined by $\chi^g
(a)=\chi(gag^{-1}), \,\  a\in H.$\\
If $g=(g_0,\alpha)\in G$ with $g_0\in SL(2,{\mathbb R})$ and
$\alpha\ in H,$ it is easy to check that for each $(\lambda,
\mu)\in {\mathbb R}^2$,
\begin{equation}
\chi^g_{(\lambda,\mu)}=\chi_{(\lambda,\mu)g_0}. \end{equation}
We
see easily from (4.2) that the $G$-orbits in ${\hat H}\cong
{\mathbb R}^2$ consist of two orbits $\Omega_0,\,\Omega_1$ given
by
$$\Omega_0=\{(0,0)\},\quad \Omega_1={\mathbb R}^2-\{(0,0)\}.$$
We observe that $\Omega_0$ is the $G$-orbit of $(0,0)$ and
$\Omega_1$ is the $G$-orbit of any element $(\lambda, \mu)\not=
0.$

\medskip

Now we choose the element $\delta = \chi_{(1, 0)}$ of ${\hat H}$.
That is, $\delta(x, y)=e^{2\pi ix}$ for all $(x,y)\in {\mathbb
R}^2.$ It is easy to check that the stabilizer of $\chi_{(0,0)}$
is $G$ and the stabilizer $G_{\delta}$ of $\delta$ is given by
$$G_{\delta}=\left\{ \left( \begin{pmatrix} 1 & 0\\c &
1\end{pmatrix},\,\alpha\right)\,\big|\ c\in {\mathbb R},\
\alpha\in {\mathbb R}^{(1,2)}\right\}.$$
We see that $H$ is
regularly embedded. This means that for every $G$-orbit $\Omega$
in ${\hat H}$ and for every $\sigma\in\Omega$ with stabilizer
$G_{\sigma}$ of $\sigma$, the canonical bijection
$G_{\sigma}\backslash G\,\longrightarrow\,\Omega$ is a
homeomorphism.

\medskip

According to G. Mackey\,([18]), we obtain
\vspace{0.1in}\\
\noindent {\bf Theorem 4.1.} {\it The irreducible unitary
representations of $G$ are the following\,:
\begin{itemize}
\item[{\rm (a)}] The irreducible unitary representations $\pi$, where the
restriction of $\pi$ to $H$ is trivial and the restriction of
$\pi$ to $SL(2, {\mathbb R})$ is an irreducible unitary
representation of $SL(2, {\mathbb R})$. For the unitary dual of
$SL(2,{\mathbb R})$, we refer to {\rm [7]} or {\rm [15], p. 123}.
\item[{\rm (b)}] The representations $\pi_{(r)}=\text{Ind}_{G_{\delta}}^G
\sigma_r\,(r\in {\mathbb R})$ induced from the unitary character
$\sigma_r$ of $G_{\delta}$ defined by $$\sigma_r\left( \left(
\begin{pmatrix} 1 & 0\\c &
1\end{pmatrix},\,(\lambda,\mu)\right)\right)=\delta(rc+\lambda)=e^{2\pi
i(rc+\lambda)},\ \ c,\lambda,\mu\in {\mathbb R}.$$ \end{itemize}}
\vspace{0.05in}
\noindent {\it Proof.} The proof of the above
theorem can be found in [22], p. 850. \hfill$\square$

\bigskip

We put $$W_1=\,\left(\,\begin{pmatrix}
0&1\\0&0\end{pmatrix},\,(0,0)\,\right),\ \ \
W_2=\,\left(\,\begin{pmatrix}
0&0\\1&0\end{pmatrix},\,(0,0)\,\right),\ \
W_3=\,\left(\,\begin{pmatrix} 1 & 0\\0 &
-1\end{pmatrix},\,(0,0)\right)$$ and
$$W_4=\,\left(\,\begin{pmatrix} 0 &0\\0&0\end{pmatrix},\,(1,0)\,\right),\
\ \ \ W_5=\,\left(\,\begin{pmatrix}
0&0\\0&0\end{pmatrix},\,(0,1)\,\right).$$ Clearly $W_1,\cdots,W_5$
form a basis of ${\mathfrak g}$.
\vspace{0.1in}\\
\noindent {\bf Lemma\,4.2.} {\it We have the following relations.
$$[W_1, W_2]=W_3,\quad [W_3, W_1]=2W_1,\quad [W_3, W_2]=-2W_2,$$
$$[W_1, W_4]=0,\quad [W_1, W_5]=-W_4,\quad [W_2, W_4]=W_5,\quad
[W_2, W_5]=0,$$ $$[W_3, W_4]=W_4,\quad [W_3, W_5]=-W_5,\quad [W_4,
W_5]=0.$$}
\vspace{0.05in}\\
\noindent {\it Proof.} The proof follows from an easy computation.
\hfill$\square$

\medskip

Let ${\mathfrak g}_{{\mathbb C}}={\mathfrak g}\otimes_{{\mathbb
R}}{\mathbb C}$ be the complexfication of ${\mathfrak g}$. We put
$${\mathfrak k}_{{\mathbb C}}={\mathbb C}\,(W_1-W_2),\quad\ \ {\mathfrak p}_{\pm}={\mathbb C}\,(W_3\pm
i(W_1+W_2)).$$ Then we have
$${\mathfrak g}_{{\mathbb C}}={\mathfrak k}_{{\mathbb C}}+{\mathfrak p}_+ +{\mathfrak p}_{-},\quad
[{\mathfrak k}_{{\mathbb C}},{\mathfrak p}_{\pm}]\subset
{\mathfrak p}_{\pm},\quad {\mathfrak p}_-={\overline{{\mathfrak
p}_+}}.$$ We note that ${\mathfrak k}_{{\mathbb C}}$ is the
complexification of the Lie algebra ${\mathfrak k}$ of $K.$

\medskip

We set $\mathfrak {a}={\mathbb R}\, W_3$. By Lemma 4.2, the roots
of ${\mathfrak g}$ relative to $\mathfrak {a}$ are given by $\pm
e,\,\pm 2e,$ where $e$ is the linear functional $e:\mathfrak
{a}\longrightarrow {\mathbb C}$ defined by $e(W_3)=1.$ The set
$\Sigma^+=\{e,\,2e\}$ is the set of positive roots of ${\mathfrak
g}$ relative to $\mathfrak {a}$. We recall that for a root
$\alpha$, the root space ${\mathfrak g}_{\alpha}$ is defined by
$${\mathfrak g}_{\alpha}=\{ X\in {\mathfrak g}\,\vert\
[H,X]=\alpha(H)X\ \text{for\ all}\ H\in \mathfrak {a}\,\}.$$ Then
we see easily that
$${\mathfrak g}_e={\mathbb R}\,W_4,\quad {\mathfrak g}_{-e}={\mathbb R}\,W_5,\quad
{\mathfrak g}_{2e}={\mathbb R}\,W_1,\quad{\mathfrak
g}_{-2e}={\mathbb R}\,W_2$$ and
$${\mathfrak g}={\mathfrak g}_{-2e}\oplus {\mathfrak g}_{-e}\oplus\mathfrak {a}\oplus
{\mathfrak g}_e\oplus{\mathfrak g}_{2e}.$$
\vspace{0.1in}\\
\noindent {\bf Proposition 4.3.} {\it The Killing form $B$ of
${\mathfrak g}$ is given by
\begin{equation}
B((X_1,Z_1),(X_2,Z_2))=\,5\,\sigma(X_1X_2), \end{equation} where
$(X_1,Z_1), (X_2,Z_2)\in {\mathfrak g}$ with $X_1,X_2\in
{\mathfrak s}{\mathfrak l}(2,{\mathbb R})$ and $Z_1,Z_2\in
{\mathbb R}^{(1,2)}.$ Hence the Killing form is highly
nondegenerate. The adjoint representation $\text{Ad}$ of $G$ is
given by
\begin{equation}
\text{Ad}((g,\alpha))(X, Z)=(gXg^{-1},\,(Z-\alpha\,^tX)\,^tg),
\end{equation}
where $(g,\alpha)\in G$ with $g\in SL(2,{\mathbb
R}),\,\alpha\in{\mathbb R}^{(1,2)}$ and $(X, Z)\in{\mathfrak g}$
with $X\in{\mathfrak s}{\mathfrak l}(2,{\mathbb R}),\,Z\in{\mathbb
R}^{(1,2)}.$}
\vspace{0.05in}\\
\noindent {\it Proof.} The proof follows immediately from a direct
computation.\hfill$\square$

\bigskip

An Iwasawa decomposition of the group $G$ is given by
\begin{equation}
G=NAK, \end{equation}
where $$N=\,\left\{\,\left(\,\begin{pmatrix} 1 & x\\
0 & 1\end{pmatrix},\, a\right)\in G\,\Big|\ x\in {\mathbb R},\
a\in {\mathbb R}^{(1,2)}\,\right\}$$ and
$$A=\,\left\{\,\left(\,\begin{pmatrix} a & 0\\ 0 & a^{-1}
\end{pmatrix},\,0\,\right)\in G\,\Big|\ a>0\ \right\}.$$
An Iwasawa decomposition of the Lie algebra ${\mathfrak g}$ of $G$
is given by $${\mathfrak g}={\mathfrak n}+ {\mathfrak
a}+{\mathfrak
k},$$ where $${\mathfrak n}=\,\left\{\,\left(\,\begin{pmatrix} 0 & x\\
0 &0
\end{pmatrix},\,Z\right)\in {\mathfrak g}\,\Big|\ x\in{\mathbb R},\ Z\in {\mathbb R}^{(1,2)}
\,\right\}$$ and $${\mathfrak a}=\,\left\{\,\left(\,\begin{pmatrix} x & 0\\
0 & -x
\end{pmatrix},\,0\right)\in {\mathfrak g}\,\Big|\ x\in{\mathbb R}\,\right\}.$$
In fact, $\mathfrak {a}$ is the Lie algebra of $A$ and ${\mathfrak
n}$ is the Lie algebra of $N$.

\medskip

Now we compute the Lie derivatives for functions on $G$
explicitly.
We define the differential operators $L_k,\,R_k\,(\,1\leq k \leq
5\,)$ on $G$ by
$$L_kf({\tilde g})=\,{ { {d}\over{dt}} \Bigg|_{t=0} } f({\tilde g}\ast \exp\,tW_k)$$
and
$$R_kf({\tilde g})=\,{ { {d}\over{dt}} \Bigg|_{t=0} } f(\exp\,tW_k\ast {\tilde g}),$$
where $f\in
C^{\infty}(G)$ and ${\tilde g}\in G.$

\medskip

By an easy calculation, we get
\begin{align*}
\exp\,tW_1\,&=\,\left(\,\begin{pmatrix} 1 & t \\ 0 &
1\end{pmatrix}\,, \,(0,\,0)\,\right),\ \ \ \
\exp\,tW_2\,=\,\left(\, \begin{pmatrix} 1 & 0 \\ t & 1
\end{pmatrix}\,,\,(\,0,\,0\,)\,\right)\\
\exp\,tW_3\,&=\,\left(\,\begin{pmatrix} e^t & 0 \\ 0 & e^{-t}
\end{pmatrix}\,,\,(\,0,\,0\,)\,\right),\ \ \ \
\exp\,tW_4\,=\,\left(\,\begin{pmatrix} 0 & 0 \\ 0 & 0
\end{pmatrix}\,,\,(\,t,\,0\,)\,\right)
\end{align*}
and $$\exp\,tW_5\,=\,\left(\,\begin{pmatrix} 0 & 0\\ 0 & 0
\end{pmatrix}\,,\,(\,0,\,t\,)\,\right).$$
Now we use the following coordinates $(g,\,\alpha)$ in $G$ given
by
\begin{equation}
g\,=\,\begin{pmatrix} 1 & x\\ 0 & 1\end{pmatrix}\, \begin{pmatrix} y^{1/2} & 0
\\ 0 & y^{-1/2}\end{pmatrix}\, \begin{pmatrix} \cos\,\theta &
\sin\,\theta\\ -\sin\,\theta & \cos\,\theta\end{pmatrix}
\end{equation}
and \begin{equation}
\alpha\,=\,(\,\alpha_1,\,\alpha_2\,),
\end{equation}
where $x,\,\alpha_1,\,\alpha_2\in {\mathbb R},\ y>0$ and $0\leq
\theta < 2\pi.$ By an easy computation, we have
\begin{align*}
L_1\,&=\, y\,\cos\,2\theta{{\partial}\over{\partial
x}}\,+\,y\,\sin\,2\theta\,{{\partial}\over{\partial y}}\,+\,
\sin^2\,\theta\,{{\partial}\over{\partial\theta}}\,-\,\alpha_2\,{{\partial}\over{\partial\alpha_1}},\\
L_2\,&=\,y\,\cos\,2\theta\,{{\partial}\over{\partial
x}}\,+\,y\,\sin\,2\theta\,{{\partial}\over{\partial y}}\,-\,
\cos^2\,\theta\,{{\partial}\over{\partial\theta}}\,-\,\alpha_1\,{{\partial}\over{\partial\alpha_2}},\\
L_3\,&=\,-\,2y\,\sin\,2\theta\,{{\partial}\over{\partial
x}}\,+\,2\,y\,\cos\,2\theta\,{{\partial}\over{\partial y}}
\,+\,\sin\,2\theta\,{{\partial}\over{\partial\theta}}\,-\,\alpha_1\,{{\partial}\over{\partial\alpha_1}}\,+\,\alpha_2\,{{\partial}\over{\partial\alpha_2}},\\
L_4\,&=\,{{\partial}\over{\partial\alpha_1}},\\ L_5\,&=\,{{\partial}\over{\partial\alpha_2}},\\ R_1\,&=\,{{\partial}\over{\partial x}},\\
R_2\,&=\,(\,y^2\,-\,x^2\,)\,{{\partial}\over{\partial x}}\,-\,2\,xy\,{{\partial}\over{\partial y}}\,-\,y\,{{\partial}\over{\partial\theta}},\\
R_3\,&=\,2\,x\,{{\partial}\over{\partial x}}\,+\,2\,y\,{{\partial}\over{\partial y}},\\
R_4\,&=\,y^{-1/2}\,\cos\,\theta\,{{\partial}\over{\partial\alpha_1}}\,+\,y^{-1/2}\,\sin\,\theta\,
{{\partial}\over{\partial\alpha_2}},\\
R_5\,&=\,-\,y^{-1/2}\,(\,x\,\cos\,\theta\,+\,y\,\sin\,\theta\,)\,
{{\partial}\over{\partial\alpha_1}}\,+\,y^{-1/2}\,(\,y\,\cos\,\theta\,-\,x\,\sin\,\theta\,)\,
{{\partial}\over{\partial\alpha_2}}.
\end{align*}

\medskip

In fact, the calculation for $L_3$ and $R_5$ can be found in [22],
p. 837-839.

\bigskip

We define the differential operators ${\mathbb L}_j\,(\,1\leq
j\leq 5\,)$ on ${\mathbb H}\times{\mathbb C}$ by
$${\mathbb L}_jf(\tau,z)=\,{d\over{dt}}\Bigg|_{t=0}f(\exp\,tW_j\circ
(\tau,z)),\ \ \ 1\leq j\leq 5,$$ where $f\in C^{\infty}({\mathbb
H}\times{\mathbb C})$. Using the coordinates $\tau\,=\,x\,+\,iy$
and $z\,=\,u\,+\,iv$ with $x,y,u,v$ real and $y>0,$ we can easily
compute the explicit formulas for ${\mathbb L}_j$'s. They are
given by
\begin{align*}
{\mathbb
L}_1&=\,(x^2-y^2)\,{{\partial}\over{\partial
x}}\,+\,2xy\,{{\partial}\over{\partial
y}}\,+\,(xu-yv)\,{{\partial}\over{\partial u}}\,+\,
(yu+xv)\,{{\partial}\over{\partial v}},\\ {\mathbb L}_2&=\,-\,{{\partial}\over{\partial x}},\\
{\mathbb L}_3&=\,-2x\,{{\partial}\over{\partial x}}\,-\,2y\,{{\partial}\over{\partial y}}\,-\,u\,{{\partial}\over{\partial u}}\,-\,v\,{{\partial}\over{\partial v}},\\
{\mathbb L}_4&=\,x\,{{\partial}\over{\partial u}}\,+\,y{{\partial}\over{\partial v}},\\
{\mathbb L}_5&=\,{{\partial}\over{\partial u}}.
\end{align*}
\vspace{0.2in}\\
\noindent{\bf 5. The decomposition of $L^2(\Gamma\backslash G)$}
\setcounter{equation}{0}
\renewcommand{\theequation}{5.\arabic{equation}}
\vspace{0.1in}\\
\indent Let $R$ be the right regular representation of $G$ on the
Hilbert space $L^2(\Gamma\backslash G)$. We set $G_1 =
SL(2,{\mathbb R}).$ Then the decomposition of $R$ is given by
\begin{equation}
L^2(\Gamma\backslash G)=L^2_{\text{disc}}(\Gamma_1\backslash
G_1)\,\bigoplus\,L^2_{\text{cont}}(\Gamma_1\backslash
G_1)\,\bigoplus\,\int_{-\infty}^{\infty}{\mathcal H}_{(r)}dr,
\end{equation}
where $L^2_{\text{disc}}(\Gamma_1\backslash G_1)$\,(resp.
\,$L^2_{\text{cont}}(\Gamma_1\backslash G_1)$) is the
discrete\,(resp.\,continuous) part of $L^2(\Gamma_1\backslash
G_1)$\,(cf.\,[14], \,[15]) and ${\mathcal H}_{(r)}$ is the
representation space of $\pi_{(r)}$\,(cf. Theorem 4.1.\,(b)).

\medskip

We recall the result of Rolf Berndt\,(cf.\,[2], \,[3], \,[4]). Let
$H_{{\mathbb R}}^{(1, 1)}$ denote the Heisenberg group which is
${\mathbb R}^3$ as a set and is equipped with the following
multiplication
$$(\lambda, \mu, \kappa) \,(\lambda', \mu', \kappa')=(\lambda+\lambda', \mu+\mu', \kappa+\kappa'+\lambda\mu'-\mu\lambda').$$
We let $G^J = SL(2,{\mathbb R})\ltimes H_{{\mathbb R}}^{(1,1)}$ be
the semidirect product of $SL(2,{\mathbb R})$ and $H_{{\mathbb
R}}^{(1,1)}$, called the Jacobi group whose multiplication law is
given by $$
(M,(\lambda,\mu,\kappa))\cdot(M',(\lambda',\mu',\kappa')) =\,
(MM',(\tilde{\lambda}+\lambda',\tilde{\mu}+ \mu',
\kappa+\kappa'+\tilde{\lambda}\mu' -\tilde{\mu}\lambda'))$$ with
$M,M'\in SL(2,{\mathbb R}),
(\lambda,\mu,\kappa),\,(\lambda',\mu',\kappa') \in H_{{\mathbb
R}}^{(1,1)}$ and $(\tilde{\lambda},\tilde{\mu})=(\lambda,\mu)M'$.
Obviously the center $Z(G^J)$ of $G^J$ is given by $\{
(0,0,\kappa)\,\vert\,\kappa\in{\mathbb R}\,\}$. We denote
$$H_{{\mathbb Z}}^{(1,1)}=\{\,(\lambda,\mu,\kappa)\in H_{{\mathbb R}}^{(1,1)}\,\vert\
\lambda,\mu,\kappa \ \text{integral}\,\}.$$ We set
$$\Gamma^J=SL(2,{\mathbb Z})\ltimes H_{{\mathbb Z}}^{(1,1)},\quad K^J=K\,\times
Z(G^J).$$ R. Berndt proved that the decomposition of the right
regular representation $R^J$ of $G^J$ in $L^2(\Gamma^J\backslash
G^J)$ is given by
\begin{equation}
L^2(\Gamma^J\backslash G^J)=\Bigg(\, \bigoplus_{m,n\in {\mathbb
Z}}{\mathcal
H}_{m,n}\,\Bigg)\,\bigoplus\,\Bigg(\,\bigoplus_{\nu=\pm{\frac12}}\int_{\scriptstyle
\text{Re}\,s=0 \atop\scriptstyle \text{Im}\,s>0}{\mathcal H}_{m,
s, \nu}\,ds\,\Bigg),\end{equation} where the ${\mathcal H}_{m, n}$
is the irreducible unitary representation isomorphic to the
discrete series $\pi_{m,k}^{\pm}$ or the principal series $\pi_{m,
s, \nu}$, and the ${\mathcal H}_{m, s, \nu}$ is the representation
space of $\pi_{m,s,\nu}$\,(cf.\,[4], \,p.\,47-48). For more detail
on the decomposition of $L^2(\Gamma^J\backslash G^J)$, we refer to
[4],\,p.\,75-103.

\medskip

Since ${\mathbb H}\,\times\,{\mathbb C}=K^J\backslash
G^J=K\backslash G,$ the space of the Hilbert space
$L^2(\Gamma\backslash ({\mathbb H}\times{\mathbb C}))$ consists of
$K^J$-fixed elements in $L^2(\Gamma^J\backslash G^J)$ or $K$-fixed
elements in $L^2(\Gamma\backslash G)$. Hence we obtain the
spectral decomposition of $L^2(\Gamma\backslash ({\mathbb
H}\times{\mathbb C}))$ for the Laplacian $\Delta$ or
$\Delta_{\alpha, \beta}$\,(cf. \,Proposition 2.4 or Remark 2.5).
\vspace{0.2in}\\
\noindent{\bf 6. Remarks on Fourier expansions of Maass-Jacobi
forms} \setcounter{equation}{0}
\renewcommand{\theequation}{6.\arabic{equation}}
\vspace{0.1in}\\
\indent  We let $f:{\mathbb H}\times{\mathbb C}\longrightarrow
{\mathbb C}$ be a Maass-Jacobi form with $\Delta f=\lambda f.$
Then $f$ satisfies the following invariance relations
\begin{equation}
f(\tau+n,\,z)\,=\,f(\tau,z)\ \ \ \text{for\ all}\ n\in
{\mathbb Z} \end{equation} and
\begin{equation}
f(\tau,\,z\,+\,n_1\tau\,+\,n_2)\,=\,f(\tau,z)\ \ \ \text{for\
all}\ n_1,\,n_2\in {\mathbb Z}.
\end{equation}
Therefore $f$ is a smooth function on ${\mathbb H}\times{\mathbb
C}$ which is periodic in $x$ and $u$ with period $1.$ So $f$ has
the following Fourier series
\begin{equation}
f(\tau,z)\,=\,\sum_{n\in {\mathbb Z}}\sum_{r\in
{\mathbb Z}}\,c_{n,r}(y,v)\, e^{2\pi i(nx+ru)}.
\end{equation}
For two fixed integers $n$ and $r$, we have to calculate the
function $c_{n, r}(y, v).$ For brevity, we put $F(y,v)=\,c_{n,
r}(y, v).$ Then $F$ satisfies the following differential equation
\begin{equation}
\left[y^2\,{{\partial^2}\over{\partial
y^2}}\,+\,(y+v^2)\,{{\partial^2}\over{\partial v^2}}\,+\,2yv\,
{{\partial^2}\over{\partial y\partial
v}}\,-\,\left\{\,(ay+bv)^2\,+\,b^2y\,+\,\lambda\,\right\}
\right]\,F\,=\,0. \end{equation}
Here $a=2\pi n$ and $b=2\pi r$
are constant. We note that the function $u(y)=y^{\frac 12}
K_{s-{\frac 12}}(2\pi \vert n\vert y)$ satisfies the differential
equation (6.4) with $\lambda=s(s-1).$ Here $K_s(z)$ is the
$K$-Bessel function defined by (2.16)\,(see Lebedev\,[16] or
Watson\,[21]). The problem is that if there exist solutions of the
differential equation (6.4), we have to find their solutions
explicitly.
\vspace{0.1in}\\
\noindent{\bf Acknowledgement.} This work started while I was
staying at Department of Mathematics, Harvard University during
the fall semester in 1996. I would like to give my hearty thanks
to Professor Don Zagier for his kind advice on this work and for
pointing out some errors in the first version. I also want to my
deep thanks to Professor Rolf Berndt for his interest on this work
and for letting me know his works.
\vspace{0.2in}\\

\footnotesize{

\end{document}
\begin{thebibliography}{99}
\bibitem{} R. Berndt, {\em Some Differential Operators in the Theory of
Jacobi Forms}, IHES/M/84/10.
\bibitem{} R. Berndt, {\em The Continuous Part of $L^2(\Gamma^J \backslash G^J)$ for
the Jacobi Group $G^J$}, Abh. Math. Sem. Univ. Hamburg, {\bf
60}(1990), 225-248.
\bibitem{} R. Berndt and S. B{\"o}cherer, {\em Jacobi
Forms and Discrete Series Representations of the Jacobi Group},
Math. Z., {\bf 204}(1990), 13-44.
\bibitem{} R. Berndt and R. Schmidt, Elements of
the Representation Theory of the Jacobi Group, Birkh{\"a}user,
{\bf 163}(1998).
\bibitem{} A. Borel and H. Jacquet, {\em Automorphic forms and automorphic
representations}, Proc. Symposia in Pure Math., {\bf XXXIII(Part
1)}(1979), 189-202.
\bibitem{} D. Bump, Automorphic Forms and Representations, Cambridge University Press,
(1997).
\bibitem{} R. W. Donley, {\em Irreducible
Representations of $SL(2, {\mathbb R})$}, Proceedings of Symposia
in Pure Mathematics on Representation Theory and Automorphic
Forms, American Math. Soc., {\bf 61}(1997), 51-59.
\bibitem{} P. R. Garabedian, Partial Differential
Equations, Wiley, New York, (1964).
\bibitem{} S. Gelbart, Automorphic forms on adele
groups, Annals of Math. Studies, Princeton Univ. Press, {\bf
83}(1975).
\bibitem{} Harish-Chandra, Automorphic forms on
semi-simple Lie groups, Notes by J.G.M. Mars, Lecture Notes in
Math., Springer-Verlag, Berlin-Heidelberg-New York, {\bf
62}(1968).
\bibitem{} S. Helgason, {\em Differential operators
on homogeneous spaces}, Acta Math., {\bf 102}(1959), 239-299.
\bibitem{} S. Helgason, Groups and geometric analysis, Academic Press,
(1984).
\bibitem{} H. Iwaniec, Introduction to the spectral
theory of automorphic forms, Biblioteca de la Revista Mathem{\'
a}tica Iberoamericana, Madrid, (1995).
\bibitem{} T. Kubota, Elementary Theory of Eisenstein
Series, John Wiley and Sons, New York, (1973).
\bibitem{} S. Lang, $SL_2({\mathbb R})$, Springer-Verlag, (1985).
\bibitem{} N. N. Lebedev, Special Functions and their Applications, Dover, New York,
(1972).
\bibitem{} H. Maass, {\em {\" U}ber eine neue Art von
nichtanalytischen automorphen Funktionen und die Bestimmung
Dirichlescher Reihen durch Funktionalgleichung}, Math. Ann., {\bf
121}(1949), 141-183.
\bibitem{} G. Mackey, {\em Unitary Representations of
Group Extensions 1}, Acta Math., {\bf 99}(1958), 265-311.
\bibitem{} A. Selberg, {\em Hamonic analysis and
discontinuous groups in weakly symmetric Riemannian spaces with
applications to Dirichlet series}, J. Indian Math. Soc., {\bf
20}(1956), 47-87.
\bibitem{} A. Terras, Harmonic analysis on symmetric spaces and
applications I, Springer-Verlag, (1985).
\bibitem{} G. N. Watson, A Treatise on the Theory of Bessel Functions, Cambridge University Press,
London, (1962).
\bibitem{} J.-H. Yang, {\em On the group $SL(2, {\mathbb R})\ltimes {\mathbb
R}^{(m, 2)}$}, J. Korean Math. Soc., {\bf 40(5)}, 831-867.


\end{thebibliography}
